\documentclass[12pt]{article}
\usepackage{amssymb,latexsym,amsmath}
\usepackage{graphics}                 
\usepackage{color}                    
\usepackage{hyperref}                 

\begin{document}

\pagestyle{plain} \headheight=5mm  \topmargin=-5mm

\title{Generalized Abel-Jacobi map on Lawson homology}
\author{ Wenchuan Hu }

\maketitle
\newtheorem{Def}{Definition}[section]
\newtheorem{Th}{Theorem}[section]
\newtheorem{Prop}{Proposition}[section]
\newtheorem{Not}{Notation}[section]
\newtheorem{Lemma}{Lemma}[section]
\newtheorem{Rem}{Remark}[section]
\newtheorem{Cor}{Corollary}[section]

\def\s{\section}
\def\ss{\subsection}

\def\d{\begin{Def}}
\def\t{\begin{Th}}
\def\p{\begin{Prop}}
\def\n{\begin{Not}}
\def\la{\begin{Lemma}}
\def\r{\begin{Rem}}
\def\c{\begin{Cor}}
\def\ee{\begin{equation}}
\def\aa{\begin{eqnarray}}
\def\ya{\begin{eqnarray*}}
\def\bd{\begin{description}}

\def\ed{\end{Def}}
\def\et{\end{Th}}
\def\epo{\end{Prop}}
\def\en{\end{Not}}
\def\el{\end{Lemma}}
\def\er{\end{Rem}}
\def\ec{\end{Cor}}
\def\eee{\end{equation}}
\def\eaa{\end{eqnarray}}
\def\ey{\end{eqnarray*}}
\def\ebd{\end{description}}

\def\nn{\nonumber}
\def\bp{{\bf Proof.}\hspace{2mm}}
\def\qe{\hfill$\Box$}
\def\lj{\langle}
\def\rj{\rangle}
\def\dd{\diamond}
\def\ox{\mbox{}}
\def\lb{\label}
\def\rel{\;{\rm rel.}\;}
\def\vp{\varepsilon}
\def\ep{\epsilon}
\def\mod{\;{\rm mod}\;}
\def\exp{{\rm exp}\;}
\def\Lie{{\rm Lie}}
\def\dim{{\rm dim}}
\def\im{{\rm im}\;}
\def\Lag{{\rm Lag}}
\def\Gr{{\rm Gr}}
\def\span{{\rm span}}
\def\Spin{{\rm Spin}}
\def\sign{{\rm sign}\;}
\def\Supp{{\rm Supp}\;}
\def\Sp{{\rm Sp}\;}
\def\ind{{\rm ind}\;}
\def\rank{{\rm rank}\;}
\def\Sg{{\Sp(2n,\C)}}
\def\Na{{\cal N}}
\def\det{{\rm det}\;}
\def\dist{{\rm dist}}
\def\deg{{\rm deg}}
\def\Tr{{\rm Tr}\;}
\def\ker{{\rm ker}\;}
\def\Vect{{\rm Vect}}
\def\H{{\bf H}}
\def\K{{\rm K}}
\def\R{{\bf R}}
\def\C{{\mathbb{C}}}
\def\Z{{\mathbb{Z}}}
\def\N{{\bf N}}
\def\F{{\bf F}}
\def\Da{{\bf D}}
\def\A{{\bf A}}
\def\La{{\bf L}}
\def\x{{\bf x}}
\def\y{{\bf y}}
\def\Ga{{\cal G}}
\def\Ha{{\cal H}}
\def\L{{\cal L}}
\def\Pa{{\cal P}}
\def\Ua{{\cal U}}
\def\E{{\rm E}}
\def\J{{\cal J}}

\def\m{{\rm m}}
\def\ch{{\rm ch}}
\def\gl{{\rm gl}}
\def\Gl{{\rm Gl}}
\def\Sp{{\rm Sp}}
\def\sf{{\rm sf}}
\def\U{{\rm U}}
\def\O{{\rm O}}
\def\F{{\rm F}}
\def\P{{\rm P}}
\def\D{{\rm D}}
\def\T{{\rm T}}
\def\Sa{{\rm S}}

\begin{abstract}
We construct an Abel-Jacobi type map on the homologically trivial
part of Lawson homology groups. It generalizes the Abel-Jacobi map
constructed by Griffiths. By using a result of H. Clemens, we give
some examples of smooth projective manifolds  with  infinite
generated Lawson homology groups $L_pH_{2p+k}(X, {\mathbb{Q}})$ when
$k>0$.
\end{abstract}

\begin{center}{\bf \tableofcontents}\end {center}

\s {Introduction}{\hskip .2 in} In this paper, all varieties are
defined over $\mathbb{C}$. Let $X$ be a smooth projective variety
with dimension $n$. Recall that the \textbf{Hodge filtration}

$$\cdots\subseteq {F}^iH^k(X,{\mathbb{C}})\subseteq {F}^{i-1}H^k(X,{\mathbb{C}})
\subseteq\cdots\subseteq {F}^0H^k(X,{\mathbb{C}})=
H^k(X,{\mathbb{C}})$$ is defined by
$$ F^qH^k(X,\mathbb{C}):=\bigoplus_{i\geq q}H^{i,k-i}(X).
$$
Note that ${F}^qH^k(X,{\mathbb{C}})$ vanishes if $q>k$.

In \cite{Griffiths}, Griffiths generalized the Jacobian varieties
and the Abel-Jacobi map on smooth algebraic curves  to higher
dimensional smooth projective varieties.

{\Def The \textbf{q-th intermediate Griffiths Jacobian} of a smooth
projective variety $X$ is defined by
\begin{eqnarray*}
J^q(X):&=&H^{2q-1}(X,\mathbb{C})/\{F^qH^{2q-1}(X,\mathbb{C})+H^{2q-1}(X,\Z) \}\\
&\cong& F^{n-q+1}H^{2n-2q+1}(X,\mathbb{C})^*/H^{2q-1}(X,\Z)^*  .
\end{eqnarray*}
}

Let ${\cal Z}_p(X)$ be the space of algebraic $p$-cycles on $X$. Set
${\cal Z}^{n-p}(X)\equiv{\cal Z}_p(X)$. There is a natural map
$$cl_q:{\cal Z}^q(X)\rightarrow H^{2q}(X,\Z)$$
called \textbf{the cycle class map}. Set
$$
{\cal Z}_{n-q}(X)_{hom}:={\cal Z}^q(X)_{hom}:=\ker cl_q
$$

\medskip
{\Def The \textbf{Abel-Jacobi map}
$$
\Phi^q: {\cal Z}^q(X)_{hom}\rightarrow J^q(X)
$$
sends $\varphi\in {\cal Z}^q(X)_{hom}$ to $\Phi_{\varphi}^q$, where
$\Phi_{\varphi}^q$ is defined by
$$
\Phi_{\varphi}^q(\omega):=\int_U\omega,\quad \omega \in
F^{n-q+1}H^{2n-2q+1}(X,\mathbb{C}).
$$
Here $\varphi=\partial U$ and $U$ is an integral current of
dimension $2n-2q+1$. }

\medskip
Now let
$$ J^{2q-1}(X)_{alg}\subseteq J^{2q-1}(X)
$$
be the largest complex subtorus of $J^{2q-1}(X)$ whose tangent space
is contained in $H^{q-1,q}(X)$. It can be proved that $\Phi^q({\cal
Z}^{q}(X)_{alg})$ is a subtorus of $J^{2q-1}(X)$ contained in
$J^{2q-1}(X)_{alg}$ (cf. \cite{Voisin1}, Corollary 12.19), where
${\cal Z}^{q}(X)_{alg}\subseteq {\cal Z}^{q}(X)$ are the subset of
codimension $q$-cycles which are algebraically equivalent to zero.

\medskip
The \textbf{Griffiths group} of codimension $q$-cycles is defined to
$$
{\rm Griff}^q(X):={\cal Z}^q(X)_{hom}/{\cal Z}^q(X)_{alg}
$$

Therefore we can define the transcendental part of the Abel-Jacobi
map
\begin{equation}\label{eq:T6eq1}
\Phi^q_{tr}: {\rm Griff}^q(X) \rightarrow
J^q(X)_{tr}:=J^{2q-1}(X)/J^{2q-1}(X)_{alg}
\end{equation}
as the factorization of $\Phi^q$.

\medskip
By using this, Griffiths showed the following:

{\Th {\rm (\cite{Griffiths})}\label{sec:ATh1.1G} Let $X\subset
{\P^4}$ be a general quintic threefold, the Griffiths group ${\rm
Griff}^2(X)$ is nontrivial, even modulo torsion.}

{\Rem \label{sec:Rem1.1} Clemens has obtained further results: Under
the same assumption as those in Theorem \ref{sec:ATh1.1G}, ${\rm
Griff}^2(X)\otimes {\mathbb{Q}}$ is an infinitely generated
${\mathbb{Q}}$-vector space \cite{Clemens}. }

\medskip
In this paper, the Griffiths' Abel-Jacobi map is generalized to the
spaces of the homologically trivial part of Lawson homology groups.

{\Def  The \textbf{Lawson homology} $L_pH_k(X)$ of $p$-cycles is
defined by
$$L_pH_k(X) := \pi_{k-2p}({\cal Z}_p(X)) \quad for\quad k\geq 2p\geq 0,$$
where ${\cal Z}_p(X)$ is provided with a natural topology (cf.
\cite{Friedlander1}, \cite{Lawson1}). For general background, the
reader is referred to \cite{Lawson2}.}

\medskip
In \cite{Friedlander-Mazur}, Friedlander and Mazur showed that there
are  natural maps, called \textbf{cycle class maps}
$$ \Phi_{p,k}:L_pH_{k}(X)\rightarrow H_{k}(X). $$
Define $$L_pH_{k}(X)_{hom}:={\rm
ker}\{\Phi_{p,k}:L_pH_{k}(X)\rightarrow H_{k}(X)\}$$ and
$$L_pH_{k}(X,\mathbb{Q}):=L_pH_{k}(X)\otimes{\mathbb{Q}}.$$

\medskip
The domain of Abel-Jacobi map can be reduced to Griffiths groups as
in (\ref{eq:T6eq1}). Similarly, our generalized Abel-Jacobi map is
defined on homologically trivial part of Lawson homology groups. As
an application, we show that the non-triviality of certain Lawson
homology group.

\medskip
The main result in this paper is the following:

{\Th \label{sec:ATh1.2} Let $X$ be a smooth projective variety.
There is a well-defined map
$$\Phi:L_pH_{2p+k}(X)_{hom}\longrightarrow \bigg\{\bigoplus_{{r>
k+1,r+s=k+1}} H^{p+r,p+s}(X)\bigg\}^*\bigg/H_{2p+k+1}(X,\Z)$$ which
generalizes  Griffiths' Abel-Jacobi map defined in [G]. Moreover,
for any $p>0$ and $k\geq 0$, we find examples of projective
manifolds $X$ for which the image of the map on
$L_pH_{k+2p}(X)_{hom}$ is infinitely generated. }

\medskip
As the application of the main result together Clemens' Theorem
(Remark \ref{sec:Rem1.1}), we obtain

 {\Th \label{sec:ATh1.3} For any $k\geq 0$, there exist a projective manifold $X$ of
 dimension $k+3$ such that $L_1H_{k+2}(X)_{hom}\otimes {\mathbb{Q}}$ is
 nontrivial, in fact, infinite dimensional over $\mathbb{Q}$.
 }

\medskip
Using the Projective Bundle Theorem proved by Friedlander and Gabber
in \cite{Friedlander-Gabber}, we have the following result:

{\Th  \label{sec:ATh1.4} For any $p>0$ and $k\geq 0$,  there exist a
smooth projective variety $X$ such that $L_pH_{k+2p}(X)_{hom}\otimes
{\mathbb{Q}}$ is an infinite dimensional vector space over
$\mathbb{Q}$.}

\medskip

In section \ref{sec:A2}, we will review the minimal background
materials about Lawson homology and point out its relation to
Griffiths groups. In section \ref{sec:A3}, we give the definition
the generalized Abel-Jacobi map. In section \ref{sec:A4}, the
non-triviality of the generalized Abel-Jacobi map is proved by using
Griffiths  and Clemens' results through examples. The construction
in our examples also shows this generalized Abel-Jacobi map really
generalizes Griffiths' result in \cite{Griffiths}.

\s{Lawson homology} {\hskip .2 in}\label{sec:A2} Let $X$ be a
projective variety of dimension $n$. Denote by ${\cal C}_p(X)$ the
space of effective algebraic $p$-cycles on $X$ and by ${\cal
Z}_p(X)$ the space of algebraic $p$-cycles on $X$. There is a
natural, compactly generated topology on ${\cal C}_p(X)$ (resp.
${\cal Z}_p(X)$) and therefore ${\cal C}_p(X)$ (resp. ${\cal
Z}_p(X)$) carries a structure of an abelian topological group
(\cite{Friedlander1}, \cite{Lawson1}).

\medskip
The \textbf{Lawson homology} $L_pH_k(X)$ of $p$-cycles is
defined by
$$L_pH_k(X) := \pi_{k-2p}({\cal Z}_p(X)) \quad for\quad k\geq 2p\geq 0.$$

It has been proved by Friedlander in \cite{Friedlander1} that
$$ L_pH_k(X)\cong \underrightarrow{\lim}\pi_k({\cal C}_p(X)_{\alpha})
$$
for all $k>0$, where the limit is taken over the connected
components of ${\cal C}_p(X)$ with respect to the action of
$\pi_0({\cal Z}_p(X))$. For a detailed discussion of this
construction and its properties we refer the reader to
\cite{Friedlander-Mazur}, \S2 and \cite{Friedlander-Lawson},\S1.

\medskip
In \cite{Friedlander-Mazur}, Friedlander and Mazur showed that there
are  natural maps, called \textbf{cycle class maps}
$$ \Phi_{p,k}:L_pH_{k}(X)\rightarrow H_{k}(X) $$
where $H_{k}(X)$ is the singular homology with the integral
coefficient.

\medskip
Define $L_pH_{k}(X)_{hom}$ to be the \textbf{homologically trivial
part of Lawson homology group} $L_pH_{k}(X)_{hom}$, i.e.,
$$L_pH_{k}(X)_{hom}:={\rm
ker}\{\Phi_{p,k}:L_pH_{k}(X)\rightarrow H_{k}(X)\}.$$

It was proved by Friedlander \cite{Friedlander1} that
$L_pH_{2p}(X)\cong {\cal Z}_p(X)/{\cal Z}_p(X)_{alg}$. Therefore we
have

\begin{equation}\label{eq:T6eq2}
L_pH_{2p}(X)_{hom}\cong {\rm Griff}_p(X),
\end{equation}
where ${\rm Griff}_p(X):={\rm Griff}^{n-p}(X)$.

For general background on Lawson homology, the reader is referred to
\cite{Lawson2}.

\s {The definition of generalized Abel-Jacobi map on
$L_pH_{2p+k}(X)_{hom}$} {\hskip .2 in} \label{sec:A3} In this
section, $X$ denotes a smooth projective algebraic manifold with
dimension $n$. Now ${\cal Z}_p(X)$ is an abelian topological group
with an identity element, the ``null" $p$-cycle.

For $[\varphi]\in L_pH_{2p+k}(X)$, we can construct an integral
$(2p+k)$-cycle $c$ in $X$ with $\Phi_{p,2p+k}([\varphi])=[c]$, where
$[c]$ is the homology class of $c$.

To see how to construct $c$ from $[\varphi]$ for the case that
$p=0$, the reader is referred to \cite{Friedlander-Lawson2}.


We will use this construction several times in the following. We
briefly review the construction here.

A class

$$
[\varphi]\in L_pH_{2p+k}(X)=\underrightarrow{\lim}\pi_k({\cal
C}_p(X))
$$
is represented by a map
$$\varphi:S^k\rightarrow {\cal C}_p(X).
$$
(For $k=0$, $[\varphi]$ is represented by  a difference of such
maps.)

We may assume $\varphi$ to be piecewise linear (PL for short) with
respect to a triangulation of ${\cal C}_p(X)\supset
\Gamma_1\supset\Gamma_2\supset\cdots$ respecting the smooth
stratified structure (\cite{Hironaka75}). Here $\Gamma_i$ is a
subcomplex for every $i>0$.

Let $\varphi$ be as above and fix $s_0\in S^k$ and $x_0\in
\Supp(\varphi(s_0))\subset X$. There exist affine coordinates
$(z_1,\cdots z_p, \zeta_1,\cdots,\zeta_{n-p})$ on $X$ with $x_0=0$
such that the projection $pr_1(z,\zeta)=z$, when restricted to
$U\times U'=\{(z,w):|z|<1 ~and~ |w|<1\}$, gives a proper (hence
finite) map $pr:\Supp(\varphi(s_0))\cap (U\times U')\rightarrow U$.
Slicing this cycle $\varphi(s_0)|_{U\times U'} $ by this projection
gives a PL map $\sigma:U\rightarrow SP^d(U')$ (with respect to a
triangulation of $SP^d(U')$) for some $d$. Furthermore, given any
such a map, we can construct a cycle in $U\times U'$. (cf.
\cite{Friedlander-Lawson2}.)  Choose a finite number of such product
neighborhood $U_{\alpha}\times U_{\alpha}'$, $\alpha=1,\cdots, K$,
so that the union of $U_{\alpha}\times U_{\alpha}'(\frac{1}{2})$
covers $\Supp(\varphi(s_0)$.  After shrinking each $U_{\alpha}'$
slightly, we can find a neighborhood $\mathcal{N}$ of $s_0$ in $S^k$
such that $pr:\Supp(\varphi(s))\cap (U_{\alpha}\times
U_{\alpha}')\rightarrow U_{\alpha}$ for all $s\in \mathcal{N}$ and
for all $\alpha$. Then $\varphi$ is PL in $\mathcal{N}$ if and only
if $\sigma:\mathcal{N}\times U_{\alpha}\rightarrow
SP^d(U_{\alpha}')$ is PL for all $\alpha$. One defines the cycle
$c(\varphi)$ in each neighborhood $\mathcal{N}\times
U_{\alpha}\times U_{\alpha}'$ by graphing this extended $\sigma$.
From the construction, the cycle $c(\varphi)$ depends only on the PL
map $\varphi$. (The argument here is from
[\cite{Friedlander-Lawson}, page 370-371].)

{\Lemma \label{sec:AL3.1} The homology class
$c_\varphi:=(pr_2)_*(c(\varphi))$ is independent of the choice of PL
map $\varphi:S^k \rightarrow {\cal C}_p(X)$ in $[\varphi]$, where
$pr_2: S^k\times X\rightarrow X$ is the projection onto the second
factor.}

\medskip
\bp Suppose that $\varphi':S^k \rightarrow {\cal C}_p(X)$ is another
PL map in $[\varphi]$. Hence, we have a continuous map $H:S^k\times
[0,1]\rightarrow {\cal C}_p(X)$ such that $H|_{S^k\times
\{0\}}=\varphi$ and $H|_{S^k\times \{1\}}=\varphi'$. Furthermore,
this map can be chosen to be PL with respect to the triangulation of
${\cal C}_p(X)$. Therefore, by the same construction as above, we
obtain that an integral current $c_H:=(pr_2)_*(c(H))$.  It is clear
$\partial(c_H)=c_{\varphi}-c_{\varphi'}$ since the push-forward
$(pr_2)_*$ commutes with the boundary map $\partial$.

\qe

\medskip
Alternatively, the restriction of $\varphi$ to the interior of each
top dimensional simplicies $\Delta^k_j$ ($1\leq j\leq N$, $N$ is the
number of top dimensional simplices) gives a map $\varphi:
\Delta^k_j\rightarrow \Gamma_{n_j}$, where $\Delta^k_j$ is the
$j$-th $k$-dimensional simplex and $n_j$ is the maximum number such
that $\Gamma_{n_{j}}$ contains the image of $\varphi|_{\Delta^k_j}$.

The piecewise linear property of $\varphi$ with respect to the
stratified structure of ${\cal C}_p(X)$ have the following property:

\begin{enumerate}
\item[$(*)$]\small \sl For each $s\in \Delta^k_j$,
$\varphi(s)=\sum_ia_i(s)V_i(s)$ with the property that $a_i(s)=a_i$
is constant in $s$ and $V_i(s)$ is irreducible.
\end{enumerate}

For each $j$, $1\leq j\leq N$, set $Z_{ \Delta^k_j}:=\sum_{i} a_i Z_{i,j}^k$, where  $Z_{i,j}^k:=\{(s,
z)\in \Delta^k_j\times X|z\in V_i(s)\}$. It is clear that $Z_{
\Delta^k_j}$ is an integral current. Therefore,

\begin{equation}\label{eq:T6eq3}
\phi^k_j:=(pr_2)_*(Z_{ \Delta^k_j})
\end{equation}
is then an integral current of real dimension $2p+k$, where
$pr_2:\Delta^k_j\times X\rightarrow X$ is the projection onto the
second factor.  Set $Z(\varphi):=\sum_{j=1}^{N}\phi^k_j$.

{\Lemma \label{sec:AL3.2} The closure of $Z(\varphi)$ is an integral
cycle in $X$.}

\medskip
\bp Since $\varphi$ is piecewise linear with respect to the
triangulation of ${\cal C}_p(X)$. The image of $\varphi$ on each
$(k-1)$-dimensional simplex $\Delta^{k-1}_i$ is in $\Gamma_{m_i}$,
where  $m_i$ is the maximum number such that $\Gamma_{m_{i}}$
contains the image of $\varphi|_{\Delta^{k-1}_i}$. Each
 $\varphi|_{\Delta^{k-1}_i}$ defines a current
$\phi^{k-1}_i:=(pr_2)_*(Z_{ \Delta^{k-1}_i})$ as in
(\ref{eq:T6eq3}). The sum

$$\sum_{i} \phi^{k-1}_i$$
is zero since, for each $\phi^{k-1}_i$, there is exactly one
$\phi^{k-1}_{i'}$ such that they have the same support but different
 orientation.

\qe

\medskip
Let $c_{\varphi}$ be the total $(k+2p)$-cycle in $X$ determined by
$[\varphi]$. We will simply use $c$ instead of $c_{\varphi}$ unless it
arises confusion.

{\Rem  $c_{\varphi}$, as current, has restricted type
$c_{\varphi}=[c_{\varphi}]_{p+k,p}+[c_{\varphi}]_{p+k-1,p+1}+\cdots+[c_{\varphi}]_{p,p+k}$.}

If $c$ is homologous to zero, we denote it by $c\stackrel{hom}{\sim}
0$, i.e., $[\varphi]\rightarrow 0$ in $H_{2p+k}(X,\Z)$ under the
natural transformation $\Phi_{p,2p+k}: L_pH_{2p+k}(X)\rightarrow
H_{2p+k}(X,\Z)$ (see, e.g., \cite{Lawson2}, p.185). This condition
translates into the fact that there exists an integral topological
$(2p+k+1)$-chain $\tilde{c}$ such that $\partial\tilde{c}=c$.

We denote by ${\rm Map} (S^k, {\cal C}_p(X))$ the set of piecewise
linear maps with respect to a triangulation of ${\cal C}_p(X)$ from
the $k$-dimensional sphere to the abelian topological monoid ${\cal
C}_p(X))$ of $p$-cycles.

Set

$${\rm Map} (S^k, {\cal C}_p(X))_{hom}\subset {\rm Map} (S^k,
{\cal C}_p(X))
$$
the subset of such maps $\varphi: S^k\rightarrow {\cal C}_p(X)$
whose total cycles $c_{\varphi}$ is homologous to zero in
$H_{2p+k}(X,\Z)$. There is a natural induced compact open topology
on the space of such maps ${\rm Map} (S^k, {\cal C}_p(X))$ (see,
e.g., Whitehead \cite{Whitehead}).

Now ${\cal Z}_p(X)$ is the group completion of the topological
monoid ${\cal C}_p(X)$ (cf. \cite{Friedlander1}, \cite{Lawson1}). In
the following, we will denote by ${\rm Map} (S^k, {\cal Z}_p(X))$
the set of piecewise linear maps with respect to a triangulation of
${\cal Z}_p(X)$ from the $k$-dimensional sphere to the abelian
topological group ${\cal Z}_p(X)$ of $p$-cycles.

Let $\varphi:S^k\rightarrow {\cal Z}_p(X)$ be a PL map which is
homotopic to zero.  Hence there exists a map
$\tilde{\varphi}:D^{k+1}\rightarrow {\cal Z}_p(X)$ such that
$\tilde{\varphi}$ is PL with respect to a triangulation of ${\cal
Z}_p(X)$ and  $\tilde{\varphi}|_{S^k}=\varphi$. Then
$\tilde{\varphi}$ determines an integral current, i.e., the total
$(k+1+2p)$-chain $\tilde{c}$ such that the boundary of $\tilde{c}$
is $c$, i.e., $\partial \tilde{c}=c$. From the definition, we have
$\varphi\in {\rm Map} (S^k, {\cal Z}_p(X))_{hom}$. Denote by ${\rm
Map} (S^k, {\cal Z}_p(X))_{0}$ the subspace of ${\rm Map} (S^k,
{\cal Z}_p(X))_{hom}$ consisting of elements $\varphi$ which are
homotopic to zero.


\ss{The generalized Abel-Jacobi map on ${\rm Map} (S^k, {\cal
Z}_p(X))_{hom}$ } {\hskip .2 in}\label{sec:3.1} In this subsection,
suppose that
$$c\stackrel{hom}{\sim} 0, ~i.e., ~[\varphi]\rightarrow 0\in
H_{2p+k}(X,\Z)$$ under the natural transformation $\Phi_{p,2p+k}:
L_pH_{2p+k}(X)\rightarrow H_{2p+k}(X,\Z)$ (see, e.g.,
\cite{Lawson2}, p.185). This condition translates into the fact that
there exists an integral topological $(2p+k+1)$-chain $\tilde{c}$
such that $\partial\tilde{c}=c$.

Consider

$$ \omega\in\bigg\{ \bigoplus_{r\geq
k+1,r+s=k+1}{\cal{E}}^{p+r,p+s}\bigg\}, \quad d\omega=0 $$ and we
define

$$ \Phi_{\varphi}(\omega)=\int_{\tilde{c}}\omega.$$
We claim:

{\Prop \label{sec:AProp3.1} $\Phi_{\varphi}$ is well-defined, i.e.,
$\Phi_{\varphi}(\omega)$, as an element in

$$\bigg\{\bigoplus_{r,s\geq
0,r+s=k+1}H^{p+r,p+s}(X)\bigg\}^*\bigg /H_{2p+k+1}(X,\Z),$$ depends
only on the cohomology class of $\omega$. Here we identify
$H_{2p+k+1}(X,\Z)$ with the image of the composition

\begin{equation}\label{eq:T6eq4}
H_{2p+k+1}(X,\Z)\stackrel{\rho}{\rightarrow} H_{2p+k+1}(X,\C) \cong
H^{2p+k+1}(X,\C)^* \stackrel{\pi}{\rightarrow}\bigg
\{\bigoplus_{r,s\geq 0,r+s=k+1}H^{p+r,p+s}(X)\bigg\}^*,
\end{equation}
where $\rho$ is the coefficient homomorphism and $\pi$ is the
projection onto the subspace. }

\medskip
\bp
We need to show
\begin{enumerate}
\item For another choice of $\omega^{\prime}\in \bigoplus_{r\geq
k+1,r+s=k+1}{\cal{E}}^{p+r,p+s}$,
$\omega-\omega^{\prime}=d\alpha$, we have
$\int_{\tilde{c}}\omega=\int_{\tilde{c}}\omega^{\prime}$.

\item If $\tilde{c}^{\prime}$  is another integral topological chain
such that $\partial\tilde{c}^{\prime}=c$, then we also have
$\int_{\tilde{c}}\omega=\int_{\tilde{c}^{\prime}}\omega$, where
$\tilde{c}^{\prime}$ is the currents determined by
$\tilde{\varphi}$.
\end{enumerate}

To show the part 1), note that we can choose $\alpha$ such that
$\omega-\omega^{\prime}=d\alpha$ for some
$\alpha$ with $\alpha^{r,s}=0$ if $r\leq k+p$ by the Hodge
decomposition theorem for differential forms on $X$.
Hence
$$\int_{\tilde{c}}\omega-\int_{\tilde{c}}\omega^{\prime}=\int_{\tilde{c}}d\alpha
=\int_{c}\alpha=0$$
by the Stokes Theorem and the reason of type. This shows that the
definition of $\Phi_{\varphi}$ is independent of the cohomology
class of

$$[\omega]\in\bigg\{
\bigoplus_{r\geq k+1,r+s=k+1}H^{p+r,p+s}(X)\bigg\}.$$

\medskip
To show the part 2), note that  $\partial( \tilde{c}-
\tilde{c}^{\prime})=0$ and hence $\tilde{c}-
\tilde{c}^{\prime}=\lambda$ is an integral topological cycle and
hence $\int_{\lambda}$ lies in the image of the composition in
(\ref{eq:T6eq4}). Hence $\int_{\lambda}$ is well-defined
independently of the choice of $\tilde{c}$ such that $\partial
\tilde{c}=c$, as an element in

$$\bigg\{\bigoplus_{r,s\geq
0,r+s=k+1}H^{p+r,p+s}(X) \bigg\}^* \bigg /H_{2p+k+1}(X,\Z).$$ Hence
we obtain a well-defined element

$$\Phi_{\varphi}\in
\bigg\{\bigoplus_{r\geq
k+1,r+s=k+1}H^{p+r,p+s}(X)\bigg\}^*\bigg/H_{2p+k+1}(X,\Z).
$$
This completes the proof of the Proposition.

\qe

\medskip
Therefore by Proposition \ref{sec:AProp3.1} we have a well-defined
homomorphism
\begin{equation}\label{eq:T6eq5}
 \Phi:{\rm Map} (S^k, {\cal Z}_p(X))_{hom}\rightarrow \bigg\{\bigoplus_{{r\geq
k+1,r+s=k+1}} H^{p+r,p+s}(X)\bigg\}^*\bigg/H_{2p+k+1}(X,\Z)
\end{equation}
given by $\Phi(\varphi)=\Phi_{\varphi}$.

\ss {The restriction of $\Phi$ on ${\rm Map} (S^k, {\cal
Z}_p(X))_{0}$}{\hskip .2 in} {\label{sec:3.2}}In this subsection, we
will study the restriction of $\Phi$ in (\ref{eq:T6eq5}) to the
subspace ${\rm Map} (S^k, {\cal Z}_p(X))_{0} \subset {\rm Map} (S^k,
{\cal Z}_p(X))_{hom}$, i.e., all PL maps from $S^k$ to ${\cal
Z}_p(X)$ which are homotopic to zero. Note that the image of $\Phi$
is in $\big\{\bigoplus_{{r\geq k+1,r+s=k+1}}
H^{p+r,p+s}(X)\big\}^*\big/H_{2p+k+1}(X,\Z)$.

Let $\varphi:S^k\rightarrow {\cal Z}_p(X)$ be an element in ${\rm
Map} (S^k, {\cal Z}_p(X))_{0}$. Denote by $c$ the total
$(2p+k)$-cycle (maybe degenerated) determined by $\varphi$. Hence
there exists a map $\tilde{\varphi}:D^{k+1}\rightarrow {\cal
Z}_p(X)$ such that $\tilde{\varphi}|_{S^k}=\varphi$ and the
associated total $(2p+k+1)$-chain $\tilde{c}$ such that the boundary
of $\tilde{c}$ is $c$, i.e., $\partial \tilde{c}=c$.

The restriction of the generalized Abel-Jacobi map $\Phi$ to the
subspace of ${\rm Map} (S^k, {\cal Z}_p(X))_{0}$ is the map

$$ \Phi_0:{\rm Map} (S^k, {\cal Z}_p(X))_{0}\rightarrow \bigg\{\bigoplus_{{r\geq
k+1,r+s=k+1}} H^{p+r,p+s}(X)\bigg\}^*\bigg/H_{2p+k+1}(X,\Z).$$

Note that

$$ \tilde{c}\in \bigg\{\bigoplus_{r,s\geq
0,r+s=k+1}{\cal{E}}_{p+r,p+s}\bigg\}$$

and

$$ c\in \bigg\{\bigoplus_{r,s\geq
0,r+s=k}{\cal{E}}_{p+r,p+s}\bigg\}.$$

Hence

$$\Phi_{\varphi}(\omega)=\int_{\tilde{c}}\omega=0$$

for

$$\omega\in \bigg\{\bigoplus_{r>k+1,r+s=k+1}{\cal{E}}^{p+r,p+s}\bigg\} \quad
with\quad d\omega=0 $$ by the reason of type. Therefore
$\Phi_{\varphi}=0$ on $$\bigoplus_{r>k+1,r+s=k+1}{H}^{p+r,p+s}(X).$$
That is to say, the image of $\Phi$ on the subspace ${\rm Map} (S^k,
{\cal Z}_p(X))_{0}$ is in
$$H^{p+k+1,p}(X)^*\bigg/\bigg\{H^{p+k+1,p}(X)^*\cap
\rho(H_{2p+k+1}(X,\Z))\bigg\}.$$

\ss {The reduction of $\Phi$ to $L_pH_{2p+k}(X)_{hom}$ }{\hskip .2
in}{\label{sec:3.3}} We reduce the domain $\Phi$ to the quotient
$${\rm Map} (S^k, {\cal Z}_p(X))_{hom}/{\rm Map} (S^k, {\cal
Z}_p(X))_{0}\cong\pi_0({\rm Map} (S^k, {\cal Z}_p(X))_{hom}.$$

 Now, if there are two PL maps $\varphi:S^k\rightarrow {\cal
Z}_p(X)$ and $\varphi^{\prime}:S^k\rightarrow {\cal Z}_p(X)$ such
that $\varphi$ is homotopic to $\varphi^{\prime}$.  Denote by $c$
(resp. $c^{\prime}$) the total $(k+2p)$-cycle determined by
$\varphi$ (resp. $\varphi^{\prime}$). For

$$ \omega\in \bigg\{\bigoplus_{r>
k+1,r+s=k+1}{\cal{E}}^{p+r,p+s}\bigg\}, \quad d\omega=0,$$ since
$c-c^{\prime}\stackrel{hom}{\sim}0$, we have
$\Phi_{\varphi-\varphi^{\prime}}\omega=\Phi_{\varphi}\omega-\Phi_{\varphi^{\prime}}\omega=0$
and

$$\Phi_{\varphi}=\Phi_{\varphi^{\prime}}\in \bigg\{\bigoplus_{{r> k+1,r+s=k+1}}
H^{p+r,p+s}(X)\bigg\}^*\bigg/H_{2p+k+1}(X,\Z)$$ by the discuss in
Section \ref{sec:3.2}.

Therefore, we have a commutative diagram
$$
\begin{array}{cccc}
{\rm Map} (S^k, {\cal Z}_p(X))_{0} & \stackrel{i}{\hookrightarrow} &
{{\rm Map} (S^k, {\cal
Z}_p(X))_{hom}}& \\
 \downarrow \Phi_0 &   & \downarrow \Phi&   \\
 H^{p+k+1,p}(X)^*/H_{2p+k+1}(X,\Z) &
\stackrel{i}{\hookrightarrow} & \big\{\bigoplus_{{r\geq
k+1,r+s=k+1}} H^{p+r,p+s}(X)\big\}^*\big/H_{2p+k+1}(X,\Z).&
\end{array}
$$

From this, we reduces $\Phi$ to a map
\begin{eqnarray}\label{eq:T6eq6}
 \Phi_{tr}:\pi_0({\rm Map} (S^k, {\cal Z}_p(X))_{hom}\rightarrow
\bigg\{\bigoplus_{{r> k+1,r+s=k+1}}
H^{p+r,p+s}(X)\bigg\}^*\bigg/H_{2p+k+1}(X,\Z)
\end{eqnarray}
given by $\Phi_{tr}(\varphi)=\Phi_{\varphi}$. Here
$/H_{2p+k+1}(X,\Z)$ means modulo the image of the composition map

$$H_{2p+k+1}(X,\Z)\stackrel{\rho}{\rightarrow} H_{2p+k+1}(X,\C)=\bigg\{\bigoplus_{{r+s=2p+k+1}}
H^{r,s}(X)\bigg\}^*\rightarrow \bigg\{\bigoplus_{{r> k+1,r+s=k+1}}
H^{p+r,p+s}(X)\bigg\}^*.$$

We complete the construction of the generalized Abel-Jacobi map on
homologically trivial part in Lawson homology

$$L_pH_{2p+k}(X)_{hom}:=
\pi_0({\rm Map} (S^k, {\cal Z}_p(X))_{hom},$$ i.e., the kernel of
the natural transformation $\Phi_{p,2p+k}:L_pH_{2p+k}(X)\rightarrow
H_{2p+k}(X,\Z)$.

{\Rem This map (\ref{eq:T6eq6}) defined above is exactly the usual
Abel-Jacobi map on Griffiths group when $k=0$ since there is a
natural isomorphism $L_pH_{2p}(X)_{hom}\cong {\rm Griff}_p(X)$ (cf.
\cite{Friedlander1}). This map $\Phi$ on $L_0H_{k}(X)_{hom}$ is
trivial since $L_0H_{k}(X)_{hom}=0$ by Dold-Thom theorem (cf.
\cite{Dold-Thom}). }

{\Rem Our generalized Abel-Jacobi map has been generalized to Lawson
homology groups by the author. The range of the more generalized
Abel-Jacobi map will be certain Deligne (co)homology. The tools used
there are ``sparks" and ``differential characters" systematically
studied by Harvey, Lawson and Zweck \cite{Harvey-Lawson-Zweck} and
\cite{Harvey-Lawson}.}

{\Rem Sometimes we also use $AJ_{X}(c)$ to denote
$\Phi_{tr}(\varphi)$, where $c$ is the cycle determined by
$\varphi$.}

{\Rem Prof. M. Walker told me the Suslin Conjecture would imply the
existence of such a generalized Abel-Jacobi map, at least for smooth
projective 4-folds with $p=1$ and $k=1$ in the equation
(\ref{eq:T6eq6}). In his recent paper, Walker has defined a morphic
Abel-Jacobi map from algebraically trivial part of $p$-cycles to
$p$-th morphic Jacobian \cite{Walker}.}

\s{The non-triviality of the generalized Abel-Jacobi map } {\hskip
.2 in} \label{sec:A4} The natural question is the existence  of
smooth projective varieties such that the generalized Abel-Jacobi
map $\Phi_{tr}$ on $L_pH_{2p+k}(X)_{hom}$ is non-trivial for both
$p>0$ and $k>0$. The following example is a family of smooth
4-dimensional projective varieties $X$ with $L_1H_{3}(X)_{hom}\neq
0$, even modulo torsion.

\medskip
{\bf Example:}  Let $E$ be a smooth elliptic curve and  $Y$ be a
smooth projective algebraic variety such that the Griffiths group of
1-cycles of $Y$ tensored with ${\mathbb{Q}}$ is nontrivial. Set
$X=E\times Y$. Let $[\omega]\in H^{4,0}(X)$ be a non zero element.
By K\"unnuth formula, we have $[\omega]=[\alpha]\wedge[\beta]$ for
some $0\neq [\alpha]\in H^{1,0}(E)$ and $0\neq [\beta]\in
H^{3,0}(Y)$.

Let $\imath:S^1\rightarrow E$ be a homeomorphism onto its image
such that $\imath (S^1)\subset E$ is not homologous to zero in
$H_1(E,\Z)$. Let $\varphi:S^1\rightarrow {\cal Z}_1(X)$ be a
continuous map given by
\begin{eqnarray}\label{eq:T6eq7}
\varphi(t)=(\imath(t),W)\in {\cal Z}_1(X),
\end{eqnarray}
where $W\in {\cal Z}_1(Y)$ a fixed element such that $W$ is
homologous to zero but $W$ is not algebraic equivalent to zero,
i.e., $W\in {\rm Griff_1(Y)}$. The existence of $W$ is the
assumption. Then there exists an integral topological chain $U$ such
that $\partial U=W$. Using the notation above, the cycle $c$
determined by $\varphi$ is $\imath(S^1)\times W$. Now
$c=\imath(S^1)\times W$ is homologous to zero in $X$. Indeed,
\begin{eqnarray}
\partial(\imath(S^1)\times U)\label{eq:T6eq8}
=\partial(\imath(S^1))\times U+\imath(S^1)\times \partial U
=\imath(S^1)\times W=c
\end{eqnarray}
 Hence $\partial\tilde{c}=\imath(S^1)\times
 U+\gamma+\partial(something)$, where $\partial\gamma=0$.
Therefore we have
$$\int_{\tilde{c}}\omega=\int_{\imath(S^1)\times U}\omega
=
\bigg(\int_{\imath(S^1)}\alpha\bigg)\cdot\bigg(\int_{U}\beta\bigg).$$

{\Prop  \label{sec:AProp4.1} Suppose  $Y$ is a smooth threefold and
$W\in {\cal Z}_1(Y)$ such that the image $AJ_Y(W)$ of $W$ under the
Griffiths' Abel-Jacobi map $AJ_Y$ is non torsion in
$H^{3,0}(Y)^*/{\rm Im}H_3(Y,\Z)$. The map $\varphi$ is given by
(\ref{eq:T6eq7}) as above. Then the map $\Phi_{tr}(\varphi)\in
H^{4,0}(X)/{\rm Im} H_4(X,\Z)$ is nontrivial, even modulo torsion.}

\medskip
\bp By K\"unneth formula, we have $H^{4,0}(E\times Y)\cong
H^{1,0}(E)\otimes H^{3,0}(Y)$ and $ H_4(E\times Y,\Z)\cong
H_4(Y,\Z)\oplus \{H_1(E,\Z)\otimes H_3(Y,\Z)
 \}\oplus\{H_2(E,\Z)\otimes H_2(Y,\Z)\}$ modulo torsion.
Let $$\pi:H_4(E\times Y,\Z)\rightarrow \{H^{4,0}(E\times Y)\}^*$$
be the natural map given by $\pi(u)(\alpha\otimes
\beta)=\int_u\alpha\wedge\beta$ for $u\in H_4(E\times Y,\Z)$ and
$\alpha \in H^{1,0}(E)$ and $\beta\in H^{3,0}(Y)$. Now $\pi(u)\neq
0$ only if $u\in H_1(E,\Z)\otimes H_3(Y,\Z)$. Hence we get
$$\{H^{4,0}(E\times Y)\}^*/{\rm Im}{H_4(E\times Y,\Z)}\cong \{H^{1,0}(E)^*\otimes
 H^{3,0}(Y)^*\}/{\rm Im}\{H_1(E,\Z)\otimes H_3(Y,\Z)
  \}.$$

Therefore, by the definition of generalized Abel-Jacobi map and
(\ref{eq:T6eq8}), we have

$$AJ_Y(\imath(S^1)\times W)(\alpha\wedge\beta)=\Phi_{tr}(\varphi)(\alpha\wedge\beta)
= \bigg(\int_{\imath(S^1)}\alpha\bigg)\cdot\bigg(\int_{U}\beta\bigg)
= \bigg(\int_{\imath(S^1)}\alpha\bigg)\cdot(AJ_Y(W)(\beta))
$$
i.e., $AJ_Y(\imath(S^1)\times W)=\int_{\imath(S^1)}\otimes AJ_Y(W)$.

Note that the map $\int_{\imath(S^1)}:H^{1,0}(E)\rightarrow
\mathbb{C}$ is in the image of the embedding
$H_1(E,\Z)\hookrightarrow H^{1,0}(E)^*$. But $AJ_Y(W)$ is a
non-torsion element in $H^{3,0}(Y)^*/{\rm Im} H_3(Y,\Z)$. Now the
conclusion of the proposition is from the following lemma.

{\Lemma \label{sec:AL4.1}Let $V_m$ and $V_n$ be two
$\mathbb{C}$-vector spaces of dimension of $m$ and $n$,
respectively. Suppose that $\Lambda_m\subset V_m$, $\Lambda_n\subset
V_n$ be two lattices, respectively. If $b\in V_n$ is a non torsion
element in $V_n$, i.e., $kb$ is not zero in $\Lambda_n$ for any
$k\in \Z^*$, then $a\otimes b$ is not in $\Lambda_m\otimes\Lambda_n$
for any $0\neq a\in \Lambda_m$. }

\medskip
\bp  Set rank$(\Lambda_m)=m_0$, rank$(\Lambda_n)=n_0$. Let
$\{e_i\}^{m_0}_{i=1}$,  $\{f_j\}^{n_0}_{j=1}$ be two integral
basis of $\Lambda_m$, $\Lambda_n$, respectively. If the conclusion
in the lemma fails, then

\begin{eqnarray} \label{eq:T6eq9} a\otimes b=\sum_{i=1}^{m_0}\sum_{j=1}^{n_0}k_{ij}e_i\otimes
 f_j,
\end{eqnarray}
for some $k_{ij}\in \Z$. By taking the conjugation, we can suppose
that $V_m$ and $V_n$ are real vector spaces with lattices
$\Lambda_m$ and $\Lambda_n$, respectively.

Suppose that $a=\sum_{i=1}^{m_0} k_ie_i$, where $k_i\in\Z$,
$i=1,\cdots, m_0$ are not all zeros. The equation (\ref{eq:T6eq9})
reads as
$$\sum_1^{m_0} k_ie_i\otimes b=\sum_{i=1}^{m_0}\sum_{j=1}^{n_0}k_{ij}e_i\otimes f_j$$
i.e.,
$$\sum_{i=1}^{m_0} e_i\otimes \big(k_ib-\sum_{j=1}^{n_0}k_{ij} f_j\big)=0.$$

Since $\{e_i\}^{m_0}_{i=1}$ is a basis in $\Lambda_m$ and hence they
are linearly independent over $\mathbb{R}$ in $V_m$, we get
\begin{eqnarray} \label{eq:T6eq10}k_ib-\sum_{j=1}^{n_0}k_{ij} f_j=0
\end{eqnarray}
for any $i=1,2,\cdots,m_0$. By assumption, at least one of $k_i$ is
nonzero since $a$ is nonzero vector in $V_m$. The equation
(\ref{eq:T6eq10}) contracts to the assumption that $kb$ is not in
$\Lambda_n$ for any $k\in \Z^*$. This completes the proof of the
lemma and hence the proof of the proposition. \qe

\qe

\medskip
More generally, we have the following Proposition

{\Prop \label{sec:AProp4.2} Suppose  $Y$ is a smooth threefold such
that the image $AJ_Y({\rm Griff_1(Y)}) $ of  ${\rm Griff_1(Y)}$
under the Griffiths' Abel-Jacobi map $AJ_Y$  tensored by
${\mathbb{Q}}$ is are infinitely dimensional ${\mathbb{Q}}$-vector
space over ${\mathbb{Q}}$ in $\{H^{3,0}(Y)^*/{\rm
Im}H_3(Y,\Z)\}\otimes { \mathbb{Q}}$. For each $W\in {\rm
Griff_1(Y)}$, The map $\varphi_{W}$ is given by (\ref{eq:T6eq7}) as
above. Then the image
$$\bigg\{\Phi_{tr}(\varphi_{W})\bigg|W\in {\rm
Griff_1(Y)} \bigg\}\otimes {\mathbb{Q}}\subset
\bigg\{H^{4,0}(X)/{\rm Im} H_4(X,\Z)\bigg\}\otimes {\mathbb{Q}}$$ is
an infinite dimensional ${\mathbb{Q}}$-vector space.}

\medskip
\bp We only need to show that:
\begin{enumerate}
\item[$(*)$]Let $N>0$ be an integer and $W_1,\cdots, W_N\in{\rm
Griff}_1(Y)$ be $N$ linearly independent elements under Griffiths
Abel-Jacobi map. Then $\varphi_{W_1},\cdots, \varphi_{W_N}\in
L_1H_3(E\times Y)_{hom}\otimes {\mathbb{Q}}$ are linearly
independent even under the generalized Abel-Jacobi map.
\end{enumerate}

The claim $(*)$ follows easily from Proposition \ref{sec:AProp4.1}
above since if $\varphi_{W_1},\cdots, \varphi_{W_N}$ are linearly
dependent implies that $W_1,\cdots, W_N$ are linearly dependent by
Proposition \ref{sec:AProp4.1}. This contradicts to the assumption.

\qe

\medskip

 Now for suitable choice of the 3-dimensional projective $Y$, for
example, the general quintic hypersurface in $\P^4$ (cf.
\cite{Griffiths}) or the Jacobian of a general algebraic curve with
genus 3 (cf. \cite{Ceresa}) and the 1-cycle $W$ whose image under
Abel-Jacobi map is nonzero, in fact, it is infinitely generated for
general quintic hypersurface in $\P^4$ (cf. \cite{Clemens}). Recall
the definition of Abel-Jacobi map, $AJ_Y(W)=\int_{U}$ module lattice
$H^3(Y,\Z)$, we have $\int_{U}\beta\neq 0$ for this choice of $W$
and some nonzero $[\beta]\in H^{3,0}(Y)$.

This example also gives an affirmative answer the following
question:

\medskip
 \noindent {\bf Question}: Can one show that
$L_pH_{2p+j}(X)_{hom}$ is nontrivial or even infinitely generated
for some projective variety $X$ where $j>0$ ?

{\Rem From the proof of the above propositions, we see that the
non-triviality of Griffiths' Abel-Jacobi map on $Y$ implies the
non-triviality of the generalized Abel-Jacobi map on homologically
trivial part of certain Lawson homology groups for $X$, i.e., all
the Abel-Jacobi invariants can be found by generalized Abel-Jacobi
map. In \cite{Clemens}, Clemens showed the for general quintic
3-fold, the image of the Griffiths group under the Griffiths'
Abel-Jacobi map is infinitely generated, even modulo torsion.}

{\Rem  Friedlander  proved in \cite{Friedlander2} the non-triviality
of $L_rH_{2p}(X)_{hom}$ for certain complete intersections by using
Nori's method in \cite{Nori}, which is totally different the
construction here.  There is no claim of any kind of infinite
generated property of Lawson homology in his paper.}

 {\Rem  Nori \cite{Nori} has generalized Theorem 0.1 and has shown that
even the Griffiths' Abel-Jacobi map is trivial on some Griffiths
group but the Griffiths group itself is nontrivial, even non
torsion. By using a total different, explicit and elementary
construction, the author has constructed  singular rational
4-dimensional projective varieties $X$ such that $L_1H_3(X)_{hom}$
is infinitely generated \cite{author}. But the Able-Jacobi map is
not defined on singular projective variety (at least I don't know).
}

\medskip
From the proof of Proposition \ref{sec:AProp4.1}, we observe that,
for $Y$ as above, and $M$ is a projective manifold, if there is a
map $i:S^k\rightarrow M$ such that
$$\int_{i(S^k)}:H^{k,0}(M)\rightarrow \mathbb{C}$$
is non-trivial as element in $\{H^{k,0}(M)\}^*$, then the value of
the generalized Abel-Jacobi map $\Phi_{tr}$ at
$\varphi:S^k\rightarrow {\cal Z}_1(X)$  defined by
$$\varphi(t)=(i(t),W)\in {\cal Z}_1(M\times Y)$$
is non-trivial, even modulo torsion.

\medskip
Note that if the {\sl complex Hurewicz homomorphism} $\rho\otimes
{\mathbb{C}}: \pi_k(X)\otimes{\mathbb{C}}\rightarrow
H_k(M,{\mathbb{C}})$ is surjective or even a little weaker
condition, i.e., the composition
$$
\pi_k(X)\otimes{\mathbb{C}}\rightarrow
H_k(M,{\mathbb{C}})\rightarrow \{H^{k,0}(M)\}^*$$
 is surjective, we have the non-triviality of
the map $\int_{i(S^k)}:H^{k,0}(M)\rightarrow \mathbb{C}$ if
$H^{k,0}(M)\neq 0$. Here the map $\pi:H_k(M,{\mathbb{C}})\rightarrow
\{H^{k,0}(M)\}^*$ is the Poincar\'{e} duality the projection
$H^k(M,{\mathbb{C}})\rightarrow H^{k,0}(M)$ in Hodge decomposition.

As a direct application to the Main Theorem in [\cite{DGMS}, \S 6]
and also Theorem 14 in \cite{Neisendorfer-Taylor}, we have the
following result on higher dimensional hypersurface.

{\Prop \label{sec:AProp4.3}Let $M$ be a smooth hypersurface in
$\P^{n+1}$ and $n>1$. Then the composition map
$$
\pi_k(X)\otimes{\mathbb{C}}\rightarrow
H_k(M,{\mathbb{C}})\rightarrow \{H^{k,0}(M)\}^*$$ is surjective for
any simply connected K\"{a}hler manifolds. }

\medskip
Therefore we obtain the following result:

{\Th \label{sec:ATh4.1}For any $k\geq 0$, there exist a projective
manifold $X$ of dimension $k+3$ such that
$L_1H_{k+2}(X)_{hom}\otimes {\mathbb{Q}}$ is \emph{nontrivial} or
even \emph{infinite} dimensional over $\mathbb{Q}$. } \qe

\medskip
By using the Projective Bundle Theorem in \cite{Friedlander-Gabber},
we get the following result:

{\Th \label{sec:ATh4.2}  For any $p>0$ and $k\geq 0$, there is a
smooth projective variety $X$ such that $L_pH_{k+2p}(X)_{hom}\otimes
{\mathbb{Q}}$ is \emph{infinite} dimensional vector space over
$\mathbb{Q}$.} \qe

\medskip
\s*{Acknowledge} I would like to express my gratitude to my advisor,
Blaine Lawson, for proposing this problem and for all his help. I
also want to thank Professor Mark Walker, for his useful remark and
interesting on this topic.

\medskip
\noindent
Department of Mathematics\\
Massachusetts Institute of Technology\\
Room 2-304\\
77 Massachusetts Avenue\\
Cambridge, MA  02139 \\
Email: wenchuan@math.mit.edu


\begin{thebibliography}{AAAA}
%
%




\bibitem[C]{Ceresa} G. Ceresa, {\sl $ C$ is not algebraically
equivalent to $C\sp{-}$ in its Jacobian.} Ann. of Math. (2) 117
(1983), no. 2, 285--291.

\bibitem[Cl]{Clemens} H. Clemens, {\sl Homological equivalence, modulo
algebraic equivalence, is not finitely generated.} Inst. Hautes
\'{e}tudes Sci. Publ. Math. No. 58, (1983), 19--38 (1984).

\bibitem[DGMS]{DGMS} P. Deligne; P. Griffiths; J. Morgan and
D. Sullivan, {\sl Real homotopy theory of K\"{a}hler manifolds.}
Invent. Math. 29 (1975), no. 3, 245--274.

\bibitem[DT]{Dold-Thom} A. Dold and R. Thom, {\sl Quasifaserungen und
unendliche symmetrische Produkte.} (German) Ann. of Math. (2) 67
1958 239--281.


\bibitem[F1]{Friedlander1} E. Friedlander, {\sl Algebraic cycles, Chow
varieties, and Lawson homology.}  Compositio Math. 77 (1991), no. 1,
55--93.

\bibitem[F2]{Friedlander2} E. Friedlander,  {\sl Relative Chow
correspondences and the Griffiths group.} (English. English, French
summary) Ann. Inst. Fourier (Grenoble) 50 (2000), no. 4, 1073--1098.

\bibitem[FL1]{Friedlander-Lawson} E. Friedlander and B. H. Lawson, {\sl A
theory of algebraic cocycles}. Ann. of Math. (2) 136 (1992), no. 2,
361--428.

\bibitem[FL2]{Friedlander-Lawson2} E. Friedlander and B. H. Lawson, {\sl
Graph mappings and poincar\'{e} duality}, preprint.

\bibitem[FG]{Friedlander-Gabber} E. Friedlander and O. Gabber, {\sl Cycle
spaces and intersection theory. Topological methods in modern
mathematics} (Stony Brook, NY, 1991), 325--370, Publish or Perish,
Houston, TX, 1993.

\bibitem[FM]{Friedlander-Mazur} E. Friedlander and B. Mazur,
{\sl Filtrations on the homology of algebraic varieties. With an
appendix by Daniel Quillen.}  Mem. Amer. Math. Soc. 110 (1994), no.
529, x+110 pp.


\bibitem[G]{Griffiths} P. Griffiths, {\sl On the periods of certain
rational integrals I, II}. Ann. of Math. (2) 90(1969), 460-495;
ibid. (2) 90(1969) 496--541.


\bibitem[HL]{Harvey-Lawson} R. Harvey;  B. Lawson, {\sl D-bar
sparks, I}.  arXiv.org:math.DG/0512247

\bibitem[HLZ]{Harvey-Lawson-Zweck} R. Harvey;  B. Lawson and J. Zweck, {\sl
The de Rham-Federer theory of differential characters and character
duality.}  Amer. J. Math. 125 (2003), no. 4, 791--847

\bibitem[Hi]{Hironaka75} H. Hironaka, Triangulations of algebraic sets, Proc. of
Symposia in Pure Math 29 (1975), 165-185.

\bibitem[H]{author} W. Hu, {\sl Infinitely generated Lawson homology
groups on some rational projective varieties}.
arXiv.org:math.AG/0602517


\bibitem[L1]{Lawson1} B. Lawson,  {\sl Algebraic cycles and homotopy
theory.} Ann. of Math. {\bf 129}(1989), 253-291.

\bibitem[L2]{Lawson2} B. Lawson,  {\sl Spaces of algebraic cycles.} pp.
137-213 in Surveys in Differential Geometry, 1995 vol.2,
International Press, 1995.


\bibitem[NT]{Neisendorfer-Taylor} J. Neisendorfer and L. Taylor, {\sl Dolbeault
homotopy theory.} Trans. Amer. Math. Soc. 245 (1978), 183--210.

\bibitem[N]{Nori} M. Nori,  {\sl Algebraic cycles and
Hodge-theoretic connectivity.} Invent. Math. 111 (1993), no. 2,
349--373.

\bibitem[V1]{Voisin1} C. Voisin, {\sl Hodge theory and complex algebraic
geometry. I.} Translated from the French original by Leila Schneps.
Cambridge Studies in Advanced Mathematics, 76. Cambridge University
Press, Cambridge, 2002. x+322 pp. ISBN 0-521-80260-1

\bibitem[Wa]{Walker} M. Walker, {\sl The morphic Abel-Jacobi map.}
\url{http://www.math.uiuc.edu/K-theory/740/}

\bibitem[W]{Whitehead} G. W. Whitehead, {\sl Elements of homotopy theory.}
Graduate Texts in Mathematics, 61. Springer-Verlag, New York-Berlin,
1978. xxi+744 pp. ISBN 0-387-90336-4

\end{thebibliography}
\end{document}